\documentclass[11pt,reqno]{amsart}

\usepackage[T1]{fontenc}
\usepackage{lmodern}
\usepackage{microtype}
\usepackage{placeins}

\usepackage{amsmath,amssymb,amsthm,mathtools}
\usepackage{booktabs}

\usepackage[
  backend=biber,
  style=alphabetic,
  maxbibnames=10,
  giveninits=true
]{biblatex}
\usepackage{xcolor}
\usepackage[
   colorlinks=true, linkcolor=blue!45!black, citecolor=blue!45!black,
  urlcolor=blue!45!black ]{hyperref}
\usepackage{orcidlink}

\numberwithin{equation}{section}

\newtheorem{theorem}{Theorem}[section]
\newtheorem{proposition}[theorem]{Proposition}
\newtheorem{lemma}[theorem]{Lemma}
\newtheorem{corollary}[theorem]{Corollary}

\theoremstyle{definition}
\newtheorem{definition}[theorem]{Definition}
\newtheorem{algorithm}[theorem]{Algorithm}
\theoremstyle{remark}
\newtheorem{remark}[theorem]{Remark}
\theoremstyle{plain}
\newtheorem{mainthm}{Theorem}

\newcommand{\R}{\mathbb R}
\newcommand{\C}{\mathbb C}
\newcommand{\Hh}{\mathbb H}
\newcommand{\Z}{\mathbb Z}
\newcommand{\Mul}{\mathsf M}
\newcommand{\Lp}{\mathsf L}
\newcommand{\Had}{\mathcal H}
\newcommand{\Ctn}{C_2^{\,n}}
\newcommand{\Rk}{\operatorname{R}}

\title{Quasilinear multiplication in Cayley--Dickson loop algebras}
\author[Harrison Lemley]{Harrison Lemley\,\orcidlink{0009-0004-5016-2152}}
\address{Department of Mathematics, Purdue University, West Lafayette, IN 47907,
USA}
\email{hlemley@purdue.edu}
\email{hlemley@algebraity.com}
\thanks{This work is part of the TANUKI-CD research program on Cayley--Dickson
algebras (\url{https://www.tanuki-cd.com}).}
\date{November~26, 2026}


\begin{document}

\begin{abstract}
Multiplication in the standard Cayley--Dickson algebra $A_n$ of rank $N := 2^n$
over a unital commutative ring $R$ with $2 \in R^\times$ admits a bilinear
algorithm requiring $O(N\log N)$ arithmetic operations in $R$ and $O(N)$
auxiliary storage. This paper uses that algorithm to develop a quasilinear
algorithm for multiplication in the loop algebra $R[Q_n]$ of the Cayley--Dickson
loop $Q_n$, reducing the $O(N^2)$ operations of the naive method to
$O(N\log N)$ operations and $O(N)$ auxiliary storage. This algorithm is also
applied to the loop algebras of central products $Q_n * Q_m$ of two
Cayley--Dickson loops with $O(nmD)$ operations and $O(D)$ storage, where
$D:=2^{n+m}$. A multiplication algorithm is also derived for finite fields of
characteristic $2$ for both loop families by lifting. Finally, an MIT-licensed
Python implementation \texttt{pycdloop} is provided with benchmarks and
floating-point error measurements.
\end{abstract}

\maketitle

\section{Introduction}

Multiplication in the standard Cayley--Dickson algebra $A_n$ of rank $N:=2^n$
over a unital commutative ring $R$ admits a bilinear algorithm requiring
$O(N\log N)$ arithmetic operations in $R$ when $2\in R^\times$
\cite{Lemley2026}, improving the previous uniform quadratic bound.
The paper in which this algorithm is
presented leaves open the task of applying this algorithm to algebras
constructed from $A_n$ which might similarly admit a quasilinear algorithm by
reducing multiplication to the product in $A_n$.

One such algebra is the loop algebra $R[Q_n]$ of the Cayley--Dickson loop
\[
Q_n=\{\pm e_p:0\le p<N\}\,,
\]
the loop of order $2N$ formed by the signed basis elements of the standard
Cayley--Dickson algebra $A_n$
\cite{Kirshtein2012Aut,Culbert2007,ChapmanEtAl2024}. Direct multiplication of
two elements in this algebra uses $4N^2$ multiplications in $R$, one for each
pair of loop elements. Group and loop algebras have been used to construct codes
and cryptosystems \cite{GonzalezEtAl2004, MarkovEtAl2020}, and a fast algorithm
for multiplication in $R[Q_n]$ may enable exploratory work towards these
applications for this algebra.

The contribution of the present paper is precisely such an algorithm. Here the
loop algebra $R[Q_n]$ decomposes as the direct sum of the group algebra
$R[\Ctn]$ and $A_n$, and the quasilinear algorithm for multiplication in $A_n$
of \cite{Lemley2026} is combined with a Walsh--Hadamard transform and
coordinatewise multiplication on the first summand to achieve an $O(N\log N)$
algorithm for multiplication in $R[Q_n]$. A similar decomposition expresses
$R[Q_n * Q_m]$ as the direct sum of $R^{2^{n+m}}$ and $A_n \otimes_R A_m$ where
$Q_n * Q_m$ is the central product of two Cayley--Dickson loops
\cite{ChapmanEtAl2024, ChapmanLevin2026}. These decompositions follow from
general results on loop algebras with a central involution
\cite{deBarros1994LoopRings} and are not claimed as a contribution of this
paper; rather, the contribution here is their organization into an efficient
algorithm for computing products in these algebras.

These decompositions further lead to corollaries regarding the bilinear
complexity of multiplication in these loop algebras. Over a field $F$ of
characteristic different from $2$, the bilinear complexity of multiplication in
$F[Q_n]$ exceeds that of $A_n$ by exactly $N$
(Proposition~\ref{prop:loop-rank}), and that of the product in $F[Q_n * Q_m]$
exceeds that of $A_n \otimes_F A_m$. The relevance of these results is that a
new tensor rank bound for either of these algebras would lead to a corresponding
bound for that of $A_n$ or $A_n \otimes_F A_m$ respectively.

The contribution of this paper culminates in and is represented by the following
theorems, which claim the existence of the quasilinear algorithms discussed
here.

\begin{samepage}
\begin{mainthm}\label{thm:main-a}
Let $R$ be a unital commutative ring with $2\in R^\times$, and for
$n\in\mathbb{Z}^{>0}$,
let $Q_n$ be the signed-basis loop of the standard Cayley--Dickson
algebra $A_n$ of rank $N$, and let $R[Q_n]$ be the corresponding loop
algebra over $R$. Then multiplication in $R[Q_n]$ can be computed using
\[
O(N\log N)
\]
arithmetic operations in $R$ and $O(N)$ auxiliary
$R$-coefficients. Algorithm~\ref{alg:loop} computes multiplication
with this time and space complexity.
\end{mainthm}
\end{samepage}

Write $N_n:=2^n$ and $N_m:=2^m$ for the ranks of the corresponding
Cayley--Dickson algebras.

\begin{samepage}
\begin{mainthm}\label{thm:main-b}
Let $R$ be a unital commutative ring with $2\in R^\times$, let
$n,m\in\mathbb{Z}^{>0}$, and let $Q_n,Q_m$ be the standard
Cayley--Dickson loops with their central product $Q_n * Q_m$.
Put $D:=N_nN_m$. Then multiplication in $R[Q_n * Q_m]$ can be computed
using
\[
O(D\log N_n\log N_m) = O(nmD)
\]
arithmetic operations in $R$ and $O(D)$ auxiliary $R$-coefficients.
Algorithm~\ref{alg:central} computes multiplication with this time
and space complexity.
\end{mainthm}
\end{samepage}

The remainder of this paper is organized as follows.
Section~\ref{sec:preliminaries} gives prerequisite and background information:
the bilinear algorithm for $A_n$, the Walsh--Hadamard transform, a general
decomposition result that applies to both $R[Q_n]$ and $R[Q_n * Q_m]$, and the
effect of scalar summands on bilinear complexity. Sections~\ref{sec:loop}
and~\ref{sec:central} prove Theorems~\ref{thm:main-a} and~\ref{thm:main-b}
respectively and give exact identities for the bilinear complexities over fields
of characteristic different from $2$. Section~\ref{sec:char2} shows that, over
fields of characteristic $2$, neither loop algebra splits as a direct sum of two
nonzero ideals; it gives a multiplication algorithm by lifting to a finite ring
of odd characteristic. Section~\ref{sec:implementation} presents the
implementation, benchmarks, and measured floating-point errors. Finally,
Sections~\ref{sec:outlook} and~\ref{sec:disclosures} give the outlook and
disclosures respectively.

\section{Preliminaries}
\label{sec:preliminaries}

\subsection{Cayley--Dickson algebras and their bilinear algorithm}

Let $R$ be a unital commutative ring with $2\in R^\times$, and for $n \in
\mathbb{Z}^{\ge0}$, let $A_n:=A_n(R)$ be the Cayley--Dickson algebra, free over
$R$ of rank $N := 2^n$, with every doubling parameter equal to $-1$, so that
$A_0:=R$ and for each $n\ge1$,
\[
A_n:=A_{n-1}\oplus A_{n-1}\,,
\]
whose elements are $(a, b) \in A_n$ for any $a, b \in A_{n-1}$, with
\begin{equation}\label{eq:CD-product}
\overline{(a,b)}:=(\bar a,-b)\,,\qquad
(a,b)(c,d):=\bigl(ac-\bar d\,b,\;da+b\bar c\bigr)\,,
\end{equation}
where $\bar r=r$ for $r\in R$ \cite[Section~2]{Lemley2026}. Elements of $A_n$
are represented with the recursively concatenated standard basis
$e_0=1,e_1,\ldots,e_{N-1}$, in which the coordinates of $(a,b)$ are those of $a$
followed by those of $b$. The doubling formula \eqref{eq:CD-product} is defined
over $\Z$ and hence defines $A_n(S)$ over every unital commutative ring $S$. The
structure constants of $A_n(R)$ in this basis are the images of those of
$A_n(\Z)$, and the latter are given by the signed-XOR rule
\cite{Bales2009,RenZhao2023}:
\begin{equation}\label{eq:signed-xor}
e_pe_q=\sigma_n(p,q)\,e_{p\oplus q}\,,\qquad\sigma_n(p,q)\in\{1,-1\}\,,
\end{equation}
where $p\oplus q$ is the bitwise exclusive or of $p$ and $q$. Write $\Ctn$ for
the elementary abelian $2$-group with integer labels $\{0,1,\ldots,N-1\}$ under
the $\oplus$ operation, where $C_2$ is the cyclic group of two elements.

\begin{theorem}[Scalar schedule {\cite{Lemley2026}}]\label{thm:schedule}
For every $n\in\mathbb{Z}^{>0}$, there is an evaluation $\mathcal E_n$
of the product of two elements $a,a'\in A_n$ in standard coordinates,
the scalar schedule of $A_n$, with the following properties.
\begin{enumerate}
\item[\textup{(i)}] The instructions of $\mathcal E_n$ are additions,
subtractions, multiplications by fixed elements of $R$, and $p_n$
multiplications of two input-dependent quantities. There are
$R$-linear forms $L_i,M_i\colon A_n\to R$ and fixed elements
$w_i\in A_n$, for $1\le i\le p_n$, such that the $i$th input-dependent
product is
\[
P_i=L_i(a)M_i(a')\,,
\]
with $L_i(a)$ computed from the first operand alone and $M_i(a')$ from
the second alone, and such that the output is formed from the products
as $\sum_iP_iw_i$. In particular,
\begin{equation}\label{eq:bilinear-identity}
aa'=\sum_{i=1}^{p_n}L_i(a)M_i(a')\,w_i\qquad(a,a'\in A_n)\,.
\end{equation}
\item[\textup{(ii)}] Writing $\ell_n$ for the number of additions,
subtractions, and multiplications by fixed elements of $R$ in
$\mathcal E_n$, and $\tau_n:=p_n+\ell_n$,
\begin{equation}\label{eq:scalar-schedule-cost}
\tau_n=O(nN)\,.
\end{equation}
\item[\textup{(iii)}] The schedule admits a sequential evaluation
using $O(N)$ auxiliary $R$-coefficients.
\end{enumerate}
\end{theorem}

\begin{proof}
Take $\mathcal E_n$ to be \cite[Algorithm~3.5]{Lemley2026} with every
doubling parameter equal to $-1$ and coefficient algebra
\[
S:=R[t]/(t^2+1)\,,
\]
which is free with basis $\{1,\iota\}$, where $\iota$ is the class of $t$
and $\iota^2=-1$. Take
\[
\rho_k:=\iota\,,\qquad \pi(a+b\iota):=a\qquad(a,b\in R)\,.
\]
Expand arithmetic in $S$ into operations in $R$ using the choices of
Remark~\ref{rem:products}. This computes the product of $A_n$ by
\cite[Proposition~3.6]{Lemley2026}. Property (i) is
\cite[Corollary~4.3]{Lemley2026} together with conditions (1)--(3) in
its proof, (ii) is \cite[Corollary~4.2]{Lemley2026}, and (iii) is
\cite[Proposition~4.5]{Lemley2026}.
\end{proof}

\begin{remark}[Number of products]\label{rem:products}
Use the level-two base case of \cite[Algorithm~3.5]{Lemley2026}.
For $a,b,c,d\in R$, multiplication in $S$ uses the three products
$ac$, $bd$, and $(a+b)(c+d)$ through the identity
\[
(a+b\iota)(c+d\iota)
=(ac-bd)+\bigl((a+b)(c+d)-ac-bd\bigr)\iota\,.
\]
For $n\ge3$, the top level uses the one-call root of
\cite[Corollary~3.4]{Lemley2026}. The two coordinate projections in the basis
$\{1,\iota\}$ recover both required combinations of lower-level alternating
products from a single recursive call over $S$. Let $c_j$ count the
input-dependent products in $R$ for one level-$j$ call to the
alternating-product routine over $S$. The level-two call uses six products in
$S$ and each higher call makes two half-size calls and uses $6(2^{j-1}-1)$
further products in $S$, so $c_2=18$ and
\[
c_j=2c_{j-1}+18(2^{j-1}-1)\qquad(j\ge3)\,.
\]
Induction gives
\[
c_j=(9j-18)2^j+18\qquad(j\ge2)\,.
\]
The root therefore uses one level-$(n-1)$ call and
$6(2^{n-1}-1)$ further products in $R$. Reconstruction uses
$3N-2$ products, so
\[
\begin{aligned}
p_n&=c_{n-1}+6(2^{n-1}-1)+3N-2\\
   &=(9n-15)2^{n-1}+10\,.
\end{aligned}
\]
The direct counts $p_1=4$ and $p_2=16$ satisfy the same formula.
\end{remark}

For the rest of this paper, write $\Mul_n(a,a')$ for the output of $\mathcal
E_n$ on $(a,a')$.

\subsection{The Walsh--Hadamard transform}

Let $\Had_N$ be the $N\times N$ matrix with entries
\[
(\Had_N)_{a,b}:=(-1)^{a\cdot b}\,,
\]
where $a\cdot b$ is the dot product of binary digits modulo $2$ for two integer
indices $a$ and $b$. Multiplication in the group algebra $R[\Ctn]$ in its group
basis
\[
\{[p]:0\le p<N\}
\]
is dyadic convolution
\[
(c*c')_k=\sum_{p\oplus q=k}c_pc'_q
\]
and expanding gives
\[
\Had_N(c*c')=(\Had_Nc)\odot(\Had_Nc')\,,
\]
where $\odot$ is the coordinate-wise product. Since $\Had_N^2=NI_N$
and $N=2^n$ is invertible in $R$, the transform is an $R$-algebra
isomorphism with
\[
\Had_N\colon R[\Ctn]\longrightarrow R^N\,,\qquad
\Had_N^{-1}=N^{-1}\Had_N\,.
\]
It is evaluated in place in $n$ stages of $N/2$ butterflies
\[
(u,v)\mapsto(u+v,u-v)\,,
\]
using $nN$ additions and subtractions \cite{FinoAlgazi1976}.

\subsection{Loops with a central involution}

A loop is a set $L$ equipped with a binary operation
$*\colon L\times L\to L$ and an identity $1\in L$ such that, writing
$ax:=a*x$, the maps
\[
x\longmapsto ax\,,\qquad x\longmapsto xa
\]
are bijective for every $a\in L$.
An element of $L$ is central if it commutes with every element of $L$
and associates with every pair of elements of $L$ in each position. For
a finite loop $L$, the loop algebra $R[L]$ is the free $R$-module with
basis $\{[x]:x\in L\}$ and the product $[x][y]:=[xy]$ for $x,y\in L$
extended bilinearly.

Both loops treated in this paper consist of the signed basis elements of an
algebra, and their loop algebras decompose by the same argument, which is given
once here. Let $L$ be a finite loop and let $z\in L$ be a central element of
order $2$. For $w\in L$ write $-w:=zw$, and for $\eta\in\{1,-1\}$ write $\eta w$
for $w$ or $-w$ when $\eta=1$ or $\eta=-1$ respectively. Since $z$ is central
and $z^2=1$,
\begin{equation}\label{eq:signs}
(\eta v)(\theta w)=\eta\theta\,(vw)
\qquad(v,w\in L,\ \eta,\theta\in\{1,-1\})\,.
\end{equation}
The subset $\{1,z\}$ is therefore a normal subloop of $L$, and $G:=L/\{1,z\}$ is
a loop whose elements are the pairs $\{w,-w\}$. Choose a representative $u_t\in
t$ for each $t\in G$, with $u_1=1$. Every element of $L$ then has a unique
expression $\eta u_t$ with $\eta\in\{1,-1\}$ and $t\in G$. Since $u_su_t$ lies
in $st$, for each $s, t$ pair there is a unique sign $\sigma(s,t)\in\{1,-1\}$
such that
\begin{equation}\label{eq:transversal}
u_su_t=\sigma(s,t)\,u_{st}
\end{equation}
for all $s,t\in G$, and $\sigma(1,t)=\sigma(t,1)=1$.

In $R[L]$, denote the basis elements $w, -w \in L$ as $[w]$ and $[-w]$
respectively for each $w, -w \in L$. Note that the sign ``$-$'' in the basis
elements is formal, such that $[-w]$ and $-[w]$ are different elements. Twisting
the product of $R[G]$ by $\sigma$ gives the free $R$-module $R^\sigma[G]$ with
basis $\{\xi_t:t\in G\}$ and the product
\[
\xi_s\xi_t:=\sigma(s,t)\,\xi_{st}
\]
extended bilinearly, whose unit is $\xi_1$. Every $x\in R[L]$ has a
unique expression
\[
x=\sum_{t\in G}\bigl(a_t[u_t]+b_t[-u_t]\bigr)\,,\qquad a_t,b_t\in R\,,
\]
and for such $x$ define
\begin{equation}\label{eq:split}
\begin{aligned}
\kappa(x)&:=\sum_{t\in G}(a_t+b_t)\,[t]\in R[G]\,,\\
\varepsilon(x)&:=\sum_{t\in G}(a_t-b_t)\,\xi_t\in R^\sigma[G]\,.
\end{aligned}
\end{equation}
Then $\kappa$ sends both $[u_t]$ and $[-u_t]$ to $[t]$, while
$\varepsilon$ sends them to $\xi_t$ and $-\xi_t$ respectively.

\begin{lemma}[Splitting]\label{lem:splitting}
Let $L$, $z$, $G$, $(u_t)_{t\in G}$, and $\sigma$ be as above, and
suppose that $2\in R^\times$. Then the map $(\kappa,\varepsilon)$ is an
$R$-algebra isomorphism
\[
R[L]\longrightarrow R[G]\oplus R^\sigma[G]\,,
\]
with the direct sum carrying componentwise operations. For
\[
c=\sum_{t\in G}c_t[t]\,,\qquad
d=\sum_{t\in G}d_t\xi_t\,,
\]
the inverse sends $(c,d)$ to
\begin{equation}\label{eq:loop-inverse}
\sum_{t\in G}\left(\frac{c_t+d_t}{2}[u_t]+\frac{c_t-d_t}{2}[-u_t]\right).
\end{equation}
\end{lemma}

\begin{proof}
On basis elements, $\kappa([\eta u_t])=[t]$ and
$\varepsilon([\eta u_t])=\eta\,\xi_t$. For $s,t\in G$ and
$\eta,\theta\in\{1,-1\}$, \eqref{eq:signs} and \eqref{eq:transversal}
give
\[
(\eta u_s)(\theta u_t)=\eta\theta\,\sigma(s,t)\,u_{st}\,,
\]
and therefore
\[
\begin{aligned}
\kappa\bigl([\eta u_s][\theta u_t]\bigr)&=[st]
=\kappa([\eta u_s])\,\kappa([\theta u_t])\,, \\
\varepsilon\bigl([\eta u_s][\theta u_t]\bigr)
&=\eta\theta\,\sigma(s,t)\,\xi_{st}
=\varepsilon([\eta u_s])\,\varepsilon([\theta u_t])\,.
\end{aligned}
\]
By bilinearity, $\kappa$ and $\varepsilon$ are multiplicative, and since $u_1=1$
they send the unit $[1]$ of $R[L]$ to the units $[1]$ and $\xi_1$ respectively.
Hence $(\kappa,\varepsilon)$ is an $R$-algebra homomorphism.

Write $\Phi$ for the map defined by \eqref{eq:loop-inverse}. Applied to
$\kappa(x)$ and $\varepsilon(x)$, it has coefficients
\[
\frac{(a_t+b_t)+(a_t-b_t)}{2}=a_t\,,\qquad
\frac{(a_t+b_t)-(a_t-b_t)}{2}=b_t\,,
\]
so $\Phi\circ(\kappa,\varepsilon)$ is the identity. Conversely, by
\eqref{eq:split}, $\kappa$ and $\varepsilon$ send the element
\eqref{eq:loop-inverse} to the arrays with coefficients
\[
\frac{c_t+d_t}{2}+\frac{c_t-d_t}{2}=c_t\,,\qquad
\frac{c_t+d_t}{2}-\frac{c_t-d_t}{2}=d_t\,,
\]
so $(\kappa,\varepsilon)\circ\Phi$ is the identity. Thus
$(\kappa,\varepsilon)$ is bijective with inverse $\Phi$.
\end{proof}

The lemma is an instance of a general result on loop rings. The central
normal subloop $\{1,z\}$ of $L$ satisfies
\[
L/\{1,z\}\cong G\,,
\]
by the definition of $G$, with the coset $\{\pm u_t\}$ corresponding to $t$, and
$\kappa$ is the $R$-linear extension of the quotient map. Hence
\cite[Proposition~3.1]{deBarros1994LoopRings} gives
\[
R[L]\cong R[G]\oplus\Delta_R(L:\{1,z\})\,,
\]
where $\Delta_R(L:\{1,z\})$ is the kernel of $\kappa$. The elements
\[
\delta_t:=([u_t]-[-u_t])/2
\]
satisfy $\kappa(\delta_t)=0$ and $\varepsilon(\delta_t)=\xi_t$, so by
Lemma~\ref{lem:splitting} they form an $R$-basis of this ideal, and
$\varepsilon$ identifies it with $R^\sigma[G]$.

\subsection{Bilinear complexity and scalar summands}
\label{subsec:bilinear}

Let $F$ be a field. For a finite-dimensional $F$-algebra $B$, not
necessarily associative or unital, write $\Rk_F(B)$ for the bilinear
complexity of its multiplication \cite[Chapter~14]{BurgisserEtAl1997}, the
least integer $r$ for which there are $F$-linear forms
$\alpha_i,\beta_i\colon B\to F$ and elements $w_i\in B$, for $1\le i\le r$,
such that
\[
bb'=\sum_{i=1}^{r}\alpha_i(b)\beta_i(b')\,w_i\qquad(b,b'\in B)\,.
\]
This is the tensor rank of the multiplication tensor of $B$ over $F$. An
isomorphism of $F$-algebras carries such an expression for one algebra to an
expression of the same length for the other, so $\Rk_F$ is invariant under
isomorphism. Moreover, the products of the terms of such expressions for two
algebras $B$ and $B'$ give
\begin{equation}\label{eq:rank-tensor}
\Rk_F(B\otimes_FB')\le\Rk_F(B)\cdot\Rk_F(B')\,,
\end{equation}
where $B\otimes_FB'$ carries the multiplication determined by
$(b\otimes b')(c\otimes c'):=(bc)\otimes(b'c')$. When $R=F$ is a field of
characteristic different from $2$, \eqref{eq:bilinear-identity} gives
\[
\Rk_F(A_n(F))\le p_n\,.
\]

Over a field of characteristic different from $2$, the splitting lemma
and the Walsh--Hadamard transform express each loop algebra considered
below as a direct sum of copies of the field and one further algebra.
The following lemma determines the bilinear complexity of such a sum.

\begin{lemma}[Scalar summands]\label{lem:scalar-summand}
Let $F$ be a field and let $B$ be a finite-dimensional $F$-algebra, not
necessarily associative or unital. Then, for every integer $k\ge0$,
\[
\Rk_F(F^k\oplus B)=k+\Rk_F(B)\,,
\]
with the direct sum carrying componentwise operations.
\end{lemma}

\begin{proof}
Adding the $k$ coordinatewise products of $F^k$ to an expression of length
$\Rk_F(B)$ for the multiplication of $B$ gives the inequality $\le$. Since
$F^k\oplus B=F\oplus(F^{k-1}\oplus B)$, it suffices by induction to prove
the inequality $\ge$ for $k=1$. Write the elements of $F\oplus B$ as
$(\gamma,b)$ with $\gamma\in F$ and $b\in B$, and suppose that
\[
(\gamma\gamma',bb')
=\sum_{i=1}^{r}\alpha_i(\gamma,b)\,\beta_i(\gamma',b')\,w_i
\]
with $F$-linear forms $\alpha_i,\beta_i$ on $F\oplus B$ and elements
$w_i\in F\oplus B$. Some $\alpha_i(1,0)$ is nonzero, since otherwise the
right-hand side vanishes at $(\gamma,b)=(\gamma',b')=(1,0)$, where the
left-hand side is $(1,0)$. After renumbering, $\alpha_r(1,0)\ne0$. Then
\[
\lambda(b):=-\frac{\alpha_r(0,b)}{\alpha_r(1,0)}
\]
is an $F$-linear form on $B$ with $\alpha_r(\lambda(b),b)=0$ for every
$b\in B$. Substituting
\[
(\gamma,b)=(\lambda(b),b)\,,\qquad
(\gamma',b')=(0,b')\,,
\]
and applying the projection $\pi_B\colon F\oplus B\to B$ gives
\[
bb'=\sum_{i=1}^{r-1}\alpha_i(\lambda(b),b)\,\beta_i(0,b')\,\pi_B(w_i)
\qquad(b,b'\in B)\,,
\]
in which each $\alpha_i(\lambda(b),b)$ is $F$-linear in $b$ and each
$\beta_i(0,b')$ is $F$-linear in $b'$. Hence $\Rk_F(B)\le r-1$.
\end{proof}

Lemma~\ref{lem:scalar-summand} is a special case of classical results on
the additivity of rank under direct sums \cite{JaJaTakche1986}, and over
$\R$ also of \cite[Theorem~2.1]{Li2026}, which proves additivity whenever
one summand is a real Hankel tensor. The proof is included because it is
short and applies over every field.

\section{Cayley--Dickson loop algebras}
\label{sec:loop}

Throughout this section, $R$ is a unital commutative ring with
$2\in R^\times$, $n\ge1$ is an integer, and $A_n=A_n(R)$ is the
standard Cayley--Dickson algebra of rank $N=2^n$.

\begin{definition}[Cayley--Dickson loops and loop algebras]
Define the Cayley--Dickson loop of order $2N=2^{n+1}$ to be
\[
Q_n:=\{e_p,-e_p:0\le p<N\}\subset A_n(\mathbb{Z})\,,
\]
with the multiplication of $A_n(\mathbb{Z})$. This loop is studied in
\cite{Kirshtein2012Aut,Culbert2007} and in
\cite[Example~15.2]{ChapmanEtAl2024}, where it is defined by the same
doubling. The Cayley--Dickson
loop algebra is the loop algebra $R[Q_n]$.
\end{definition}

The decomposition of Section~\ref{sec:preliminaries} is now applied
to the loop algebra $R[Q_n]$.

\begin{theorem}[Decomposition]\label{thm:decomposition}
For
\[
x=\sum_p\bigl(a_p[e_p]+b_p[-e_p]\bigr)\in R[Q_n]\,,
\]
define
\begin{equation}\label{eq:split-n}
\begin{aligned}
\kappa_n(x)&:=\sum_p(a_p+b_p)\,[p]\in R[\Ctn]\,, \\
\varepsilon_n(x)&:=\sum_p(a_p-b_p)\,e_p\in A_n(R)\,.
\end{aligned}
\end{equation}
With the direct sums carrying componentwise operations, the resulting
$R$-algebra isomorphism for each $n$ is
\[
\begin{aligned}
R[Q_n]&\xrightarrow{(\kappa_n,\varepsilon_n)}R[\Ctn]\oplus A_n(R)\\
      &\cong R^N\oplus A_n(R)\,.
\end{aligned}
\]
\end{theorem}

\begin{proof}
The element $-1=-e_0$ is central of order $2$ in $Q_n$, since $e_0$ is
the unit of $A_n(\mathbb{Z})$ and its multiplication is bilinear. By
\eqref{eq:signed-xor}, the quotient-loop isomorphism is
\[
Q_n/\{1,-1\}\longrightarrow \Ctn\,,\qquad
\{e_p,-e_p\}\longmapsto p\,.
\]
Taking $u_p:=e_p$ gives $u_0=1$ and $\sigma=\sigma_n$, so
Lemma~\ref{lem:splitting} applies. By \eqref{eq:signed-xor} again,
the $R$-algebra isomorphism
\[
R^{\sigma_n}[\Ctn]\longrightarrow A_n(R)\,,\qquad
\xi_p\longmapsto e_p
\]
identifies the second component with $A_n(R)$, giving the stated map.
Its inverse is \eqref{eq:loop-inverse} with $u_p=e_p$. The last
statement follows since $\mathcal H_N$ is an
algebra isomorphism $R[\Ctn]\to R^N$.
\end{proof}

For $n\le2$, the loop $Q_n$ is a group, and over fields of
characteristic different from $2$ the decomposition is the classical
decomposition of its group algebra; for instance,
$\R[Q_2]\cong\R^4\oplus\Hh$ for the quaternion group $Q_2$. The case
$n=4$ over such fields also appears in
\cite[Proposition~6.6(i)]{deBarrosJuriaans1996Flexible}.

Theorem~\ref{thm:decomposition} reduces multiplication in $R[Q_n]$ to
separate multiplications in its two summands: the first by the Walsh--Hadamard
transform, and the second by $\Mul_n$. In Algorithm~\ref{alg:loop},
the coordinate identifications are
\[
\begin{aligned}
(c_p)_p&\longmapsto\sum_p c_p[p]\in R[\Ctn]\,,\\
(d_p)_p&\longmapsto\sum_p d_p e_p\in A_n(R)\,.
\end{aligned}
\]

\noindent\begin{minipage}{\linewidth}
\begin{algorithm}[Loop algebra multiplication]\label{alg:loop}
For $x,y\in R[Q_n]$, compute $\Lp_n(x,y)$ by
\begin{tabbing}
\qquad\=\qquad\=\kill
\>write $x=\sum_p\bigl(a_p[e_p]+b_p[-e_p]\bigr)$ and\\
\>  $y=\sum_p\bigl(a_p'[e_p]+b_p'[-e_p]\bigr)$;\\
\>$c\leftarrow(a_p+b_p)_p$,\\
\>$d\leftarrow(a_p-b_p)_p$,\\
\>$c'\leftarrow(a_p'+b_p')_p$,\\
\>$d'\leftarrow(a_p'-b_p')_p$;\\
\>$f\leftarrow N^{-1}\Had_N\bigl((\Had_Nc)\odot(\Had_Nc')\bigr)$;\\
\>$g\leftarrow\Mul_n(d,d')$;\\
\>\textbf{return} $\sum_p\Bigl(\dfrac{f_p+g_p}{2}[e_p]
  +\dfrac{f_p-g_p}{2}[-e_p]\Bigr)$.
\end{tabbing}
\end{algorithm}
\end{minipage}\par

\begin{proposition}[Correctness and cost]\label{prop:loop-correctness}
For all $x,y\in R[Q_n]$, $\Lp_n(x,y)=xy$. Algorithm~\ref{alg:loop} uses
at most
\[
\tau_n+3nN+10N
\]
arithmetic operations in $R$, which is $O(N\log N)$, and admits a
sequential evaluation
using $O(N)$ auxiliary $R$-coefficients. Each of the $N+p_n$
multiplications of two input-dependent quantities is the product of an
$R$-linear form in the coordinates of $x$ with an $R$-linear form in
the coordinates of $y$.
\end{proposition}

\begin{proof}
By \eqref{eq:split-n},
\[
\begin{aligned}
c&=\kappa_n(x),& d&=\varepsilon_n(x),\\
c'&=\kappa_n(y),& d'&=\varepsilon_n(y)\,;
\end{aligned}
\]
then
\[
f=\kappa_n(x)\kappa_n(y)\,,
\]
since $\Had_N$ is an algebra isomorphism onto $R^N$, and
\[
g=\varepsilon_n(x)\varepsilon_n(y)
\]
by Theorem~\ref{thm:schedule}. As $\kappa_n$ and $\varepsilon_n$
are homomorphisms,
\[
(f,g)=(\kappa_n(xy),\varepsilon_n(xy))\,,
\]
and the algorithm returns
\[
(\kappa_n,\varepsilon_n)^{-1}(f,g)=xy
\]
by Theorem~\ref{thm:decomposition}.

For the cost, forming the arrays $c,d,c',d'$ and the returned coefficients uses
$6N$ additions and subtractions and $2N$ halvings, and forming $f$ via the fast
Walsh--Hadamard transform uses $3nN$ additions and subtractions, $N$
multiplications, and $N$ multiplications by $N^{-1}$. Forming $g$ uses
$\tau_n=O(nN)$ operations by Theorem~\ref{thm:schedule}(ii). The arrays
$c,d,c',d',f,g$ have length $N$, the transforms are evaluated in place, and
$\Mul_n$ admits a sequential evaluation with $O(N)$ auxiliary coefficients by
Theorem~\ref{thm:schedule}(iii). The multiplications of two input-dependent
quantities are the $N$ coordinatewise products, whose factors are coordinates of
$\Had_Nc$ and $\Had_Nc'$, and the $p_n$ products of $\mathcal E_n$, whose
factors are $L_i(d)$ and $M_i(d')$; these factors are $R$-linear in $x$ and in
$y$ respectively.
\end{proof}

By Remark~\ref{rem:products}, the number of multiplications of two
input-dependent quantities in $R$ used by Algorithm~\ref{alg:loop}
is exactly
\[
N+p_n = N+(9n-15)2^{n-1}+10\,.
\]

\begin{proof}[Proof of Theorem~\ref{thm:main-a}]
The decomposition is Theorem~\ref{thm:decomposition}, and the
correctness, arithmetic and storage bounds, and bilinearity of the
algorithm are Proposition~\ref{prop:loop-correctness}.
\end{proof}

Over a field $F$ of characteristic different from $2$, no bilinear algorithm for
$F[Q_n]$ can save a multiplication by mixing the two summands of
Theorem~\ref{thm:decomposition}.

\begin{proposition}[Bilinear complexity]\label{prop:loop-rank}
Let $F$ be a field of characteristic different from $2$. Then
\[
\Rk_F(F[Q_n])=N+\Rk_F(A_n(F))\,.
\]
\end{proposition}

\begin{proof}
By Theorem~\ref{thm:decomposition}, $F[Q_n]\cong F^N\oplus A_n(F)$ as
$F$-algebras, with the direct sum carrying componentwise operations. Since
$\Rk_F$ is invariant under isomorphism, the assertion is
Lemma~\ref{lem:scalar-summand} with $k=N$ and $B=A_n(F)$.
\end{proof}

With $R=F$, Algorithm~\ref{alg:loop} uses $N+p_n$ input-dependent products, so
by Proposition~\ref{prop:loop-rank} it exceeds the minimum by exactly
$p_n-\Rk_F(A_n(F))$. Any bilinear algorithm for $F[Q_n]$ using fewer than
$N+p_n$ products would therefore imply the existence of one for $A_n(F)$ using
fewer than $p_n$.

\begin{corollary}[Explicit bounds]\label{cor:loop-rank-bounds}
Let $F$ be a field of characteristic different from $2$.
\begin{enumerate}
\item[\textup{(i)}] For $n\ge4$,
\[
\begin{aligned}
\Rk_F(A_n(F))&\le(36n-69)2^{n-3}+10\,,\\
\Rk_F(F[Q_n])&\le N+(36n-69)2^{n-3}+10\,.
\end{aligned}
\]
\item[\textup{(ii)}] If $-1$ is a square in $F$, then for $n\ge2$,
\[
\begin{aligned}
\Rk_F(A_n(F))&\le(3n-3)2^n+4\,,\\
\Rk_F(F[Q_n])&\le N+(3n-3)2^n+4\,.
\end{aligned}
\]
\end{enumerate}
\end{corollary}

\begin{proof}
By Proposition~\ref{prop:loop-rank}, it suffices to prove the bounds for
$A_n(F)$, whose doubling parameters are all equal to $-1$. If $-1$ is a
square in $F$, then every doubling parameter is a square in $F$, and (ii)
is \cite[Corollary~4.4]{Lemley2026}. In this case (i) follows from (ii),
since
\[
(36n-69)2^{n-3}+10-\bigl((3n-3)2^n+4\bigr)=(12n-45)2^{n-3}+6>0
\qquad(n\ge4)\,.
\]
If $-1$ is not a square in $F$, then $K:=F[t]/(t^2+1)$ is a quadratic
extension of $F$ in which the class $\iota$ of $t$ satisfies
$\iota^2=-1$. Since every doubling parameter equals $-1$, the negative
of each is $1$, a square in $K$, and $\iota$ is a square root of each,
so (i) is \cite[Theorem~A.5]{Lemley2026}.
\end{proof}

The upper bound for $\Rk_F(A_n(F))$ in \textup{(i)} uses the octonion
base case of \cite[Appendix~A]{Lemley2026}. For $n\ge4$, the level-two
schedule $\mathcal E_n$ of Theorem~\ref{thm:schedule} uses
$9\cdot2^{n-3}$ more input-dependent products than this bound.

\begin{remark}[Small cases]\label{rem:loop-rank-small}
Over $\R$, the algebras $A_1(\R)\cong\C$ and $A_2(\R)\cong\Hh$ have
multiplication tensor ranks $3$ \cite[Chapter~17]{BurgisserEtAl1997} and $8$
\cite{HowellLafon1975,Lysikov2012}, respectively.
Proposition~\ref{prop:loop-rank} gives
\[
\begin{aligned}
\Rk_\R(\R[Q_1])&=2+\Rk_\R(\C)=2+3=5\,,\\
\Rk_\R(\R[Q_2])&=4+\Rk_\R(\Hh)=4+8=12\,.
\end{aligned}
\]
The first value can be confirmed independently: $Q_1$ is cyclic of
order $4$, so
\[
\R[Q_1]\cong\R[t]/(t^4-1)\,,
\]
and since $t^4-1$ has three distinct irreducible factors over $\R$, Winograd's
rank formula \cite{Winograd1977}
\[
\Rk_F(F[t]/(f))=2\deg f-r
\]
for $f$ with $r$ distinct irreducible
factors over a sufficiently large field gives $2\cdot4-3=5$.
\end{remark}

\section{Central-product loop algebras}
\label{sec:central}

Throughout this section, $R$ is a unital commutative ring with
$2\in R^\times$, $n,m\in\mathbb{Z}^{>0}$, and $A_n,A_m$ and $Q_n,Q_m$
are as in Sections~\ref{sec:preliminaries} and~\ref{sec:loop}. Write
$N_n:=2^n$, $N_m:=2^m$, and
$D:=N_nN_m$. This section mirrors the pattern of
Section~\ref{sec:loop}: the central product $Q_n*Q_m$ is realized as
the loop of signed basis elements of a tensor product, its
loop algebra is decomposed by Lemma~\ref{lem:splitting}, a
multiplication routine is given for the
one summand not already treated, and the pieces are combined into an
algorithm.

\subsection{The central product}

The tensor product $A_n\otimes_RA_m$ is an $R$-algebra under the
multiplication determined by
\[
(a\otimes b)(a'\otimes b'):=(aa')\otimes(bb')\,.
\]
Like its factors, it is in general neither commutative nor associative.
It is free of rank $D$ with standard basis
\[
e_{p,q}:=e_p\otimes e_q
\qquad(0\le p<N_n,\ 0\le q<N_m)\,,
\]
and \eqref{eq:signed-xor} gives
\begin{equation}\label{eq:central-tensor-product}
e_{p,q}e_{r,s}=\sigma_n(p,r)\sigma_m(q,s)\,e_{p\oplus r,q\oplus s}\,.
\end{equation}

\begin{definition}[Central products and their loop algebras]
\label{def:central-product}
The central product of $Q_n$ and $Q_m$ as in \cite{ChapmanLevin2026} is the set
\[
Q_n*Q_m:=\{e_{p,q},-e_{p,q}:0\le p<N_n,\ 0\le q<N_m\}
\subset A_n(\mathbb Z)\otimes_{\mathbb Z}A_m(\mathbb Z)\,,
\]
with the multiplication of $A_n(\mathbb Z)\otimes_{\mathbb Z}A_m(\mathbb Z)$.
Its loop algebra is $R[Q_n*Q_m]$.
\end{definition}

The map $(x,y)\mapsto x\otimes y$ from the direct product $Q_n\times Q_m$
onto $Q_n*Q_m$ preserves products, and
\[
x\otimes y=x'\otimes y'
\quad\Longleftrightarrow\quad
(x',y')\in\{(x,y),(-x,-y)\}\,.
\]
Since $\{(1,1),(-1,-1)\}$ is a central, hence normal, subloop of
$Q_n\times Q_m$, it follows that
\[
Q_n*Q_m\cong (Q_n\times Q_m)/\{(1,1),(-1,-1)\}\,.
\]
This loop has order $2D$, and the quotient identifies the central
involutions of the two factors. When $n,m\ge2$, this is the central
product of \cite[Section~2]{ChapmanLevin2026} over the common center
$\{\pm1\}$, with all doubling parameters equal to $-1$. When either
factor is $Q_1$, which is cyclic of order $4$, the quotient uses its
central subgroup $\{\pm1\}$. The interchange $(x,y)\mapsto(y,x)$
preserves the diagonal sign subgroup and hence induces an isomorphism
\[
Q_n*Q_m\cong Q_m*Q_n\,.
\]

\subsection{The loop algebra and its decomposition}

For $0\le p<N_n$ and $0\le q<N_m$, put
\[
h(p,q):=pN_m+q\,.
\]
Since $N_m=2^m$, the binary digits of $h(p,q)$ are those of $p$
followed by those of $q$, so $h$ is a group isomorphism
$\Ctn\times C_2^{\,m}\to C_2^{\,n+m}$:
\[
h(p\oplus r,q\oplus s)=h(p,q)\oplus h(r,s)\,.
\]
Arrays $(c_{p,q})_{p,q}$ of $D$ coefficients are stored in the order
given by $h$, so that they consist of $N_n$ consecutive rows of length
$N_m$. An array $(d_{p,q})_{p,q}$ is read as
\[
\sum_{p,q}d_{p,q}e_{p,q}\in A_n\otimes_RA_m\,.
\]
For
\[
x=\sum_{p,q}\bigl(a_{p,q}[e_{p,q}]+b_{p,q}[-e_{p,q}]\bigr)
\in R[Q_n*Q_m]\,,
\]
define
\begin{equation}\label{eq:central-split}
\begin{aligned}
\kappa_{n,m}(x)
&:=\sum_{p,q}(a_{p,q}+b_{p,q})\,[h(p,q)]\in R[C_2^{\,n+m}]\,,\\
\varepsilon_{n,m}(x)
&:=\sum_{p,q}(a_{p,q}-b_{p,q})\,e_{p,q}\in A_n\otimes_RA_m\,,
\end{aligned}
\end{equation}
and call $\varepsilon_{n,m}(x)$ the signed component of $x$.

\begin{theorem}[Decomposition]\label{thm:central-decomposition}
With the direct sums carrying componentwise operations, the resulting
$R$-algebra isomorphisms are
\[
\begin{aligned}
R[Q_n*Q_m]
&\xrightarrow{(\kappa_{n,m},\varepsilon_{n,m})}
  R[C_2^{\,n+m}]\oplus(A_n\otimes_RA_m)\\
&\cong R^D\oplus(A_n\otimes_RA_m)\,.
\end{aligned}
\]
\end{theorem}

\begin{proof}
The element $-1=-e_{0,0}$ is central of order $2$ in $Q_n*Q_m$, since
$e_{0,0}$ is the unit of $A_n(\mathbb Z)\otimes_{\mathbb Z}A_m(\mathbb Z)$
and its multiplication is bilinear. By \eqref{eq:central-tensor-product}
and the homomorphism property of $h$, the quotient-loop isomorphism is
\[
\begin{aligned}
(Q_n*Q_m)/\{1,-1\}&\longrightarrow C_2^{\,n+m}\,,\\
\{e_{p,q},-e_{p,q}\}&\longmapsto h(p,q)\,.
\end{aligned}
\]
Taking $u_{h(p,q)}:=e_{p,q}$ gives $u_0=1$ and
\[
\sigma\bigl(h(p,q),h(r,s)\bigr)=\sigma_n(p,r)\sigma_m(q,s)\,.
\]
Lemma~\ref{lem:splitting} therefore applies. By
\eqref{eq:central-tensor-product} again, the $R$-algebra isomorphism
\[
\begin{aligned}
R^\sigma[C_2^{\,n+m}]&\longrightarrow A_n\otimes_RA_m\,,\\
\xi_{h(p,q)}&\longmapsto e_{p,q}
\end{aligned}
\]
identifies the second component with $A_n\otimes_RA_m$. Under these
identifications, $\kappa$ and $\varepsilon$ become $\kappa_{n,m}$ and
$\varepsilon_{n,m}$, giving the stated map. Its inverse is
\eqref{eq:loop-inverse} with $u_{h(p,q)}=e_{p,q}$. The last statement
follows since $\Had_D$ is an algebra isomorphism
$R[C_2^{\,n+m}]\to R^D$.
\end{proof}

\subsection{Multiplication in \texorpdfstring{$A_n\otimes_RA_m$}{An tensor Am}}

Following Algorithm~\ref{alg:loop}, a product in $R[Q_n*Q_m]$ will be
computed componentwise through Theorem~\ref{thm:central-decomposition}.
The first component is handled by the Walsh--Hadamard transform as in
Section~\ref{sec:loop}; the second requires a multiplication procedure
for $A_n\otimes_RA_m$.

For an array $U=(U_{p,q})_{p,q}$, its $p$th row and the corresponding
tensor satisfy
\[
\begin{aligned}
U_p&:=\sum_{q<N_m} U_{p,q}e_q\in A_m\,,\\
U&=\sum_{p<N_n} e_p\otimes U_p\,.
\end{aligned}
\]
For an array $V$ with rows $V_r$ defined likewise,
\eqref{eq:signed-xor} gives
\[
UV=\sum_{p,r}\sigma_n(p,r)\,e_{p\oplus r}\otimes(U_pV_r)\,.
\]
Thus a product in $A_n\otimes_RA_m$ is a product in $A_n$ with
coefficients in $A_m$, and evaluating this formula directly requires
$N_n^2$ products in $A_m$. The routine below instead evaluates
the scalar schedule $\mathcal E_n$ with coefficients in $A_m$. This is
the standard extension of a bilinear algorithm to a tensor product
\cite[Chapter~14]{BurgisserEtAl1997}; since $A_m$ is in general
neither commutative nor associative, it is recorded here in the form
needed, in which each product keeps the factor coming from the first
operand on the left.

\begin{lemma}[Coefficient extension]\label{lem:coefficient-extension}
Let $B$ be an $R$-module with an $R$-bilinear multiplication, not
necessarily commutative or associative, and give $A_n\otimes_RB$ the
multiplication determined by
$(a\otimes b)(a'\otimes b'):=(aa')\otimes(bb')$. Let $L_i,M_i,w_i$, for
$1\le i\le p_n$, be as in Theorem~\ref{thm:schedule}\textup{(i)}, and
define the $R$-linear maps
$\widetilde L_i,\widetilde M_i\colon A_n\otimes_R B\to B$ by
\[
\begin{aligned}
\widetilde L_i(a\otimes b)&:=L_i(a)b\,,\\
\widetilde M_i(a\otimes b)&:=M_i(a)b\,.
\end{aligned}
\]
Then
\[
UV=\sum_{i=1}^{p_n}w_i\otimes
\bigl(\widetilde L_i(U)\,\widetilde M_i(V)\bigr)
\qquad(U,V\in A_n\otimes_RB)\,.
\]
\end{lemma}

\begin{proof}
The maps $\widetilde L_i,\widetilde M_i$ are well defined because
their defining expressions are $R$-bilinear in $(a,b)$. Both sides of
the claimed identity are $R$-bilinear in $(U,V)$, so it suffices to
take $U=a\otimes b$ and $V=a'\otimes b'$. By $R$-bilinearity of the
multiplication of $B$ and \eqref{eq:bilinear-identity},
\[
\begin{aligned}
\sum_iw_i\otimes\bigl((L_i(a)b)(M_i(a')b')\bigr)
&=\Bigl(\sum_iL_i(a)M_i(a')w_i\Bigr)\otimes(bb')\\
&=(aa')\otimes(bb')=UV.
\end{aligned}
\qedhere
\]
\end{proof}

\noindent\begin{minipage}{\linewidth}
\begin{algorithm}[Signed-component multiplication]\label{alg:signed}
For arrays $U,V$ of size $N_n$ by $N_m$, compute
$\mathsf P_{n,m}(U,V)$ by
\begin{tabbing}
\qquad\=\qquad\=\kill
\>replace each register of $\mathcal E_n$ by a\\
\>\>register of $N_m$ coefficients, read as an element of $A_m$;\\
\>load the rows $(U_p)_p$ and $(V_p)_p$ as the operands;\\
\>follow $\mathcal E_n$, performing\\
\>\>additions, subtractions, and\\
\>\>\quad multiplications by fixed elements of $R$ coordinatewise,\\
\>\>each product $P_i$ by evaluating
  $\mathcal E_m$ on its two factors,\\
\>\>\quad in the order of Theorem~\ref{thm:schedule}\textup{(i)};\\
\>\textbf{return} the output registers as the rows of an array.
\end{tabbing}
\end{algorithm}
\end{minipage}\par

\begin{lemma}[Correctness of the signed routine]
\label{lem:signed-correctness}
For all $U,V\in A_n\otimes_RA_m$ in standard coordinates,
$\mathsf P_{n,m}(U,V)=UV$.
\end{lemma}

\begin{proof}
Let $\widetilde L_i,\widetilde M_i$ be as in
Lemma~\ref{lem:coefficient-extension} with $B=A_m$. The linear
instructions of $\mathcal E_n$ act on $A_m$-valued registers by the
same $R$-linear combinations as on scalar registers, and
\[
\begin{aligned}
\widetilde L_i(U)
&=\widetilde L_i\Bigl(\sum_p e_p\otimes U_p\Bigr)
  =\sum_p L_i(e_p)U_p\,,\\
\widetilde M_i(V)
&=\widetilde M_i\Bigl(\sum_r e_r\otimes V_r\Bigr)
  =\sum_r M_i(e_r)V_r\,.
\end{aligned}
\]
Hence, by Theorem~\ref{thm:schedule}(i), Algorithm~\ref{alg:signed}
forms these factors in this order and computes their product in
$A_m$ by $\mathcal E_m$. By Lemma~\ref{lem:coefficient-extension},
the routine returns
\[
\sum_{i=1}^{p_n}w_i\otimes
\bigl(\widetilde L_i(U)\widetilde M_i(V)\bigr)=UV\,.
\qedhere
\]
\end{proof}

\subsection{The central-product algorithm}

Write $\mathsf K_{n,m}$ for the following procedure. Arrays are
stored in the order given by $h$, so that $\Had_D$ acts on them as on
vectors of length $D$.

\noindent\begin{minipage}{\linewidth}
\begin{algorithm}[Central-product multiplication]\label{alg:central}
For $x,y\in R[Q_n*Q_m]$, compute $\mathsf K_{n,m}(x,y)$ by
\begin{tabbing}
\qquad\=\qquad\=\kill
\>write $x=\sum_{p,q}\bigl(a_{p,q}[e_{p,q}]+b_{p,q}[-e_{p,q}]\bigr)$
  and\\
\>  $y=\sum_{p,q}\bigl(a'_{p,q}[e_{p,q}]+b'_{p,q}[-e_{p,q}]\bigr)$;\\
\>$c\leftarrow(a_{p,q}+b_{p,q})_{p,q}$,\\
\>$d\leftarrow(a_{p,q}-b_{p,q})_{p,q}$,\\
\>$c'\leftarrow(a'_{p,q}+b'_{p,q})_{p,q}$,\\
\>$d'\leftarrow(a'_{p,q}-b'_{p,q})_{p,q}$;\\
\>$f\leftarrow D^{-1}\Had_D\bigl((\Had_Dc)\odot(\Had_Dc')\bigr)$;\\
\>$g\leftarrow\mathsf P_{n,m}(d,d')$;\\
\>\textbf{return} $\sum_{p,q}\Bigl(\dfrac{f_{p,q}+g_{p,q}}{2}[e_{p,q}]
  +\dfrac{f_{p,q}-g_{p,q}}{2}[-e_{p,q}]\Bigr)$.
\end{tabbing}
\end{algorithm}
\end{minipage}\par

\begin{proposition}[Correctness and cost]
\label{prop:central-correctness}
For all $x,y\in R[Q_n*Q_m]$, 
\[
\mathsf K_{n,m}(x,y)=xy
\]
and Algorithm~\ref{alg:central} uses at most
\[
N_m\ell_n+p_n\tau_m+3(n+m)D+10D
\]
arithmetic operations in $R$, which is $O(D\log N_n\log N_m)$, and
admits a sequential evaluation using $O(D)$ auxiliary
$R$-coefficients. Each of its $D+p_np_m$ multiplications of two
input-dependent quantities is the product of an $R$-linear form in the
coordinates of $x$ with an $R$-linear form in the coordinates of $y$.
\end{proposition}

\begin{proof}
By \eqref{eq:central-split},
\[
\begin{aligned}
c&=\kappa_{n,m}(x),& d&=\varepsilon_{n,m}(x),\\
c'&=\kappa_{n,m}(y),& d'&=\varepsilon_{n,m}(y)\,;
\end{aligned}
\]
then
\[
f=\kappa_{n,m}(x)\kappa_{n,m}(y)\,,
\]
since $\Had_D$ is an algebra isomorphism onto $R^D$, and
\[
g=\varepsilon_{n,m}(x)\varepsilon_{n,m}(y)
\]
by Lemma~\ref{lem:signed-correctness}. As $\kappa_{n,m}$ and
$\varepsilon_{n,m}$ are homomorphisms,
\[
(f,g)=(\kappa_{n,m}(xy),\varepsilon_{n,m}(xy))\,,
\]
and the algorithm returns
\[
(\kappa_{n,m},\varepsilon_{n,m})^{-1}(f,g)=xy
\]
by Theorem~\ref{thm:central-decomposition}.

For the cost, forming $c,d,c',d'$ and the returned coefficients uses
$6D$ additions and subtractions and $2D$ halvings, and forming $f$
uses $3(n+m)D$ additions and subtractions, $D$ multiplications, and
$D$ multiplications by $D^{-1}$. In Algorithm~\ref{alg:signed}, each
of the $\ell_n$ additions, subtractions, and fixed-scalar
multiplications of $\mathcal E_n$ acts on $N_m$ coordinates, and each of
its $p_n$ input-dependent products is one evaluation of
$\mathcal E_m$ and uses $\tau_m$ operations. This gives the stated bound.
By \eqref{eq:scalar-schedule-cost}, $\ell_n,p_n\le\tau_n=O(nN_n)$
and $\tau_m=O(mN_m)$, so the bound is
\[
\begin{aligned}
O\bigl(nD+nmD+(n+m)D\bigr)
&=O(nmD)\\
&=O(D\log N_n\log N_m)\,.
\end{aligned}
\]

For the storage, the arrays $c,d,c',d',f,g$ have length $D$, and the
transforms are evaluated in place. In Algorithm~\ref{alg:signed}, the
$O(N_n)$ auxiliary registers of $\mathcal E_n$ given by
Theorem~\ref{thm:schedule}(iii) each hold at most $N_m$
coefficients, for $O(D)$ coefficients in all, and only one evaluation
of $\mathcal E_m$ is in progress at any time, using $O(N_m)$ auxiliary
coefficients by Theorem~\ref{thm:schedule}(iii).

The multiplications of two input-dependent quantities are the $D$
coordinatewise products and, within each of the $p_n$ evaluations of
$\mathcal E_m$, its $p_m$ products. The factors of the latter are
$R$-linear forms in $\widetilde L_i(d)$ and in $\widetilde M_i(d')$
respectively by Theorem~\ref{thm:schedule}(i) applied to
$\mathcal E_m$, hence $R$-linear in $x$ and in $y$.
\end{proof}

By Remark~\ref{rem:products}, Algorithm~\ref{alg:central} uses
\[
D+\bigl((9n-15)2^{n-1}+10\bigr)\bigl((9m-15)2^{m-1}+10\bigr)
\]
multiplications of two input-dependent quantities in $R$.

\begin{proof}[Proof of Theorem~\ref{thm:main-b}]
The decomposition is Theorem~\ref{thm:central-decomposition}, and the
correctness, arithmetic and storage bounds, and bilinearity of the
algorithm are Proposition~\ref{prop:central-correctness}.
\end{proof}

As in Section~\ref{sec:loop}, the first summand of
Theorem~\ref{thm:central-decomposition} contributes exactly its rank to
the bilinear complexity over a field.

\begin{proposition}[Bilinear complexity]\label{prop:central-rank}
Let $F$ be a field of characteristic different from $2$. Then
\[
\begin{aligned}
\Rk_F(F[Q_n*Q_m])&=D+\Rk_F\bigl(A_n(F)\otimes_FA_m(F)\bigr) \\
&\le D+\Rk_F(A_n(F))\,\Rk_F(A_m(F))\,.
\end{aligned}
\]
In particular, for $n,m\ge4$,
\[
\Rk_F(F[Q_n*Q_m])\le
D+\bigl((36n-69)2^{n-3}+10\bigr)\bigl((36m-69)2^{m-3}+10\bigr)\,.
\]
\end{proposition}

\begin{proof}
By Theorem~\ref{thm:central-decomposition},
$F[Q_n*Q_m]\cong F^D\oplus(A_n(F)\otimes_FA_m(F))$ as $F$-algebras, so the
equality is Lemma~\ref{lem:scalar-summand} with $k=D$. The inequality is
\eqref{eq:rank-tensor}, and the last assertion follows from
Corollary~\ref{cor:loop-rank-bounds}\textup{(i)}.
\end{proof}

Thus the minimum number of input-dependent products among bilinear algorithms
for $F[Q_n*Q_m]$ exceeds the minimum for $A_n(F)\otimes_FA_m(F)$ by exactly $D$.

\begin{remark}[More than two factors]\label{rem:k-fold}
The argument of this section is not specific to two factors. For
$n_1,\ldots,n_k\in\mathbb{Z}^{>0}$ and $D:=2^{n_1+\cdots+n_k}$, the
central product $Q_{n_1}*\cdots*Q_{n_k}$ is the loop of signed basis
elements of $A_{n_1}\otimes_R\cdots\otimes_RA_{n_k}$
\cite{ChapmanLevin2026}, so Lemma~\ref{lem:splitting} gives
\[
R[Q_{n_1}*\cdots*Q_{n_k}]\cong
R^D\oplus(A_{n_1}\otimes_R\cdots\otimes_RA_{n_k})\,.
\]
Applying Lemma~\ref{lem:coefficient-extension} with
$B=A_{n_2}\otimes_R\cdots\otimes_RA_{n_k}$ and repeating for the
remaining factors multiplies the second summand using
$O(D\,n_1\cdots n_k)$ arithmetic operations in $R$ and
$p_{n_1}\cdots p_{n_k}$ multiplications of two input-dependent
quantities, with the implied constant depending on $k$.

Over a field $F$ of characteristic different from $2$,
Lemma~\ref{lem:scalar-summand} and \eqref{eq:rank-tensor} likewise give
\[
\begin{aligned}
&\Rk_F\bigl(F[Q_{n_1}*\cdots*Q_{n_k}]\bigr)\\
&\qquad=D+\Rk_F\bigl(A_{n_1}(F)\otimes_F\cdots\otimes_F A_{n_k}(F)\bigr)\\
&\qquad\le D+\prod_{j=1}^{k}\Rk_F(A_{n_j}(F))\,.
\end{aligned}
\]
\end{remark}

\section{Characteristic 2}
\label{sec:char2}

Theorems~\ref{thm:main-a} and~\ref{thm:main-b} require $2$ to be invertible.
This section first shows that no decomposition of the kind in
Lemma~\ref{lem:splitting} is available in characteristic $2$, and then recovers
multiplication over finite fields of characteristic $2$ by lifting the
coefficients to a ring in which $2$ is invertible, multiplying there by
Algorithm~\ref{alg:loop} or~\ref{alg:central}, and reducing the result.
Throughout this section, $L$ is $Q_n$ or $Q_n*Q_m$, of order $2D$ with $D=N$ in
the first case and $D=N_nN_m$ in the second.

\subsection{The absence of a splitting}

For a field $K$ and
\[
x=\sum_{g\in L}x_g[g]\in K[L]\,,
\]
define
\[
\operatorname{aug}(x):=\sum_{g\in L}x_g\,.
\]
The map $\operatorname{aug}\colon K[L]\to K$ is a homomorphism of
$K$-algebras. As for loops in Section~\ref{sec:preliminaries}, an
element of $K[L]$ is central if it commutes with every element of
$K[L]$ and associates with every pair of elements of $K[L]$ in each
position.

\begin{proposition}[No splitting in characteristic $2$]
\label{prop:char2-no-splitting}
Let $K$ be a field of characteristic $2$ and let $x\in K[L]$. Then
$x^2$ is central in $K[L]$, the powers of $x$ do not depend on
bracketing, and
\[
x^4=\operatorname{aug}(x)^4\,.
\]
In particular, the only idempotents of $K[L]$ are $0$ and $1$, so
$K[L]$ is not a direct sum of two nonzero ideals.
\end{proposition}

\begin{proof}
Write $u_t$ for the basis elements $e_p$ of $A_n$ or $e_{p,q}$ of
$A_n\otimes A_m$, with $t=p$ or $t=h(p,q)$ respectively. By
\eqref{eq:signed-xor} and \eqref{eq:central-tensor-product}, every
element of $L$ is $u_t$ or $-u_t$ for a unique index $t$, and
$u_su_t=\pm u_{s\oplus t}$. Hence, for all $g_1,g_2,g_3\in L$,
\[
g_2g_1=\pm g_1g_2\,,\qquad
g_1(g_2g_3)=\pm(g_1g_2)g_3\,,\qquad
g_1^2=\pm1\,,
\]
since in each case the elements compared have the same index. Put
$\zeta:=[1]+[-1]$. Since $-1$ is central in $L$, the element $[-1]$ is
central in $K[L]$, and thus so is $\zeta$. Moreover,
\[
\begin{aligned}
\zeta[g]&=[g]+[-g]=\zeta[-g]\qquad(g\in L)\,,\\
\zeta^2&=2[1]+2[-1]=0\,.
\end{aligned}
\]

For every $w\in K[L]$, the element $\zeta w$ is central. Indeed, by
linearity it suffices to take $w=[g_1]$ and to compare products with
basis elements $[g_2]$ and $[g_3]$. By the centrality of $\zeta$ and
the relation $g_2g_1=\pm g_1g_2$,
\[
\begin{aligned}
(\zeta[g_1])[g_2]
&=\zeta[g_1g_2]\\
&=\zeta[g_2g_1]\\
&=[g_2](\zeta[g_1])\,.
\end{aligned}
\]
Likewise, the two bracketings of a product of $\zeta[g_1]$, $[g_2]$,
and $[g_3]$, with $\zeta[g_1]$ in any of the three positions, are
$\zeta[a]$ and $\zeta[b]$ for the two bracketings $a,b$ of the
corresponding product in $L$, and $b=\pm a$.

Now write
\[
x^2=\sum_{g\in L}x_g^2\,[g^2]
+\sum_{g,g'\in L}x_gx_{g'}\bigl([gg']+[g'g]\bigr)\,,
\]
where the second sum is over the unordered pairs of distinct elements
of $L$. The first sum is $\alpha[1]+\beta[-1]$ with
\[
\alpha+\beta=\sum_{g\in L}x_g^2=\operatorname{aug}(x)^2\,.
\]
For the mixed terms,
\[
[gg']+[g'g]=
\begin{cases}
2[gg']=0, & g'g=gg',\\
[gg']+[-(gg')]=\zeta[gg'], & g'g=-(gg').
\end{cases}
\]
Thus
\[
x^2=\alpha[1]+\beta[-1]+\zeta w
\]
for some $w\in K[L]$, a sum of three central elements, so $c:=x^2$ is
central. Every bracketing of $j\ge1$ factors equal to $x$ has value
\[
\begin{cases}
c^r,  & j=2r,\\
c^r x,& j=2r+1.
\end{cases}
\]
For $j=1$ the formula gives $x$. For $j\ge2$, each bracketing is a
product of two shorter ones. By induction, and since the central
elements of $K[L]$ form a commutative and associative subalgebra, the
powers of $c$ can be collected, and when both lengths are odd the two
remaining factors $x$ multiply to $c$. Hence the powers of $x$ do not
depend on bracketing. Since the three summands of $c$ are central,
$[-1]^2=[1]$, and $(\zeta w)^2=\zeta^2w^2=0$,
\[
x^4=c^2=\alpha^2+\beta^2=(\alpha+\beta)^2=\operatorname{aug}(x)^4\,.
\]
Finally, if $x^2=x$, then $x=x^4=\operatorname{aug}(x)^4$ lies in $K$
and is an idempotent of $K$, hence $0$ or $1$. A decomposition of
$K[L]$ as a direct sum of two nonzero ideals would give an idempotent
other than $0$ and $1$.
\end{proof}

\subsection{Multiplication by lifting}

Over a finite field of characteristic $2$, the product in $K[L]$ is
instead recovered from a product in odd characteristic. Let
$k\in\mathbb{Z}^{>0}$, let $f\in\mathbb F_2[t]$ be irreducible of degree
$k$, and put $K:=\mathbb F_2[t]/(f)$. For $x\in K[L]$, let
$\tilde x\in\Z[t][L]$ be the element each of whose coefficients is the
polynomial of degree less than $k$ with coefficients in $\{0,1\}$ that
represents the corresponding coefficient of $x$.

\begin{lemma}[Lifting]\label{lem:lifting}
For all $x,y\in K[L]$, the following hold.
\begin{enumerate}
\item[\textup{(i)}] The product $xy$ is the image of $\tilde x\tilde y$
under the map $\Z[t][L]\to K[L]$ that reduces coefficients modulo $2$
and $f$.
\item[\textup{(ii)}] Every coefficient of $\tilde x\tilde y$ is a
polynomial of degree at most $2k-2$ whose coefficients are integers
between $0$ and $2kD$.
\item[\textup{(iii)}] Let $p>2kD$ be a prime, let
$\varphi\in\mathbb F_p[t]$ be monic of degree $2k-1$, and put
$E:=\mathbb F_p[t]/(\varphi)$. Then each coefficient of
$\tilde x\tilde y$ is the polynomial of degree less than $2k-1$ with
coefficients in $\{0,1,\ldots,p-1\}$ that represents the corresponding
coefficient of the product in $E[L]$ of the images of $\tilde x$ and
$\tilde y$.
\end{enumerate}
\end{lemma}

\begin{proof}
A homomorphism of commutative rings $S\to S'$ induces a ring
homomorphism $S[L]\to S'[L]$ acting on coefficients, since the product
of two basis elements of a loop ring is a basis element. For (i), the
reduction $\Z[t]\to K$ sends each coefficient of $\tilde x$ to the
corresponding coefficient of $x$, so the induced homomorphism sends
$\tilde x$ to $x$, $\tilde y$ to $y$, and hence $\tilde x\tilde y$ to
$xy$. For (ii), the coefficient of $[g]$ in $\tilde x\tilde y$ is
\[
(\tilde x\tilde y)_g
=\sum_{\substack{g_1,g_2\in L\\g_1g_2=g}}
  \tilde x_{g_1}\tilde y_{g_2}\,.
\]
Since $L$ is a loop, there is exactly one such pair for each
$g_1\in L$, so this sum has $|L|=2D$ terms. Each term is a product of
two polynomials of degree less than $k$ with coefficients in
$\{0,1\}$, hence a polynomial of degree at most $2k-2$ in which the
coefficient of $t^j$ counts the pairs of exponents with sum $j$ at
which both factors have coefficient $1$, and so lies between $0$ and
$k$. For (iii), the reduction $\Z[t]\to E$ induces a homomorphism
sending $\tilde x\tilde y$ to the product in $E[L]$ of the images of
$\tilde x$ and $\tilde y$. By (ii), and since $2k-2<\deg\varphi$ and
$2kD<p$, each coefficient of $\tilde x\tilde y$ is a polynomial of
degree less than $2k-1$ with coefficients in $\{0,1,\ldots,p-1\}$, and
is therefore the representative of its image in $E$ described in (iii).
\end{proof}

Lemma~\ref{lem:lifting} reduces a product in $K[L]$ to one product in
$E[L]$, which Algorithm~\ref{alg:loop} or~\ref{alg:central} computes
whenever $p$ is odd.

\noindent\begin{minipage}{\linewidth}
\begin{algorithm}[Lifted multiplication]\label{alg:lift}
Fix a prime $p>2kD$ and a monic $\varphi\in\mathbb F_p[t]$ of degree
$2k-1$, and put $E:=\mathbb F_p[t]/(\varphi)$. For $x,y\in K[L]$,
compute $\mathsf T_L(x,y)$ by
\begin{tabbing}
\qquad\=\qquad\=\kill
\>$\hat x,\hat y\leftarrow$ the images of $\tilde x,\tilde y$ in $E[L]$;\\
\>$\hat z\leftarrow\Lp_n(\hat x,\hat y)$ if $L=Q_n$,\\
\>\>or $\hat z\leftarrow\mathsf K_{n,m}(\hat x,\hat y)$ if $L=Q_n*Q_m$;\\
\>$z\leftarrow$ the element of $\Z[t][L]$ whose coefficients are\\
\>\>the representatives of those of $\hat z$ described in\\
\>\>Lemma~\ref{lem:lifting}\textup{(iii)};\\
\>\textbf{return} the reduction of $z$ modulo $2$ and $f$.
\end{tabbing}
\end{algorithm}
\end{minipage}\par

\begin{proposition}[Correctness and cost]\label{prop:lift-correctness}
For all $x,y\in K[L]$, $\mathsf T_L(x,y)=xy$. Algorithm~\ref{alg:lift}
uses the arithmetic operations in $E$ of one call to
Algorithm~\ref{alg:loop} or~\ref{alg:central}, which number
$O(N\log N)$ for $L=Q_n$ and $O(D\log N_n\log N_m)$ for $L=Q_n*Q_m$,
together with $2D$ reductions modulo $2$ and $f$ of polynomials of
degree at most $2k-2$.
\end{proposition}

\begin{proof}
Since $2kD\ge4$, the prime $p$ is odd, so $E$ is a unital commutative
ring in which $2$ is invertible.
Proposition~\ref{prop:loop-correctness}
or~\ref{prop:central-correctness} with $R=E$ therefore gives
$\hat z=\hat x\hat y$, the image of $\tilde x\tilde y$ in $E[L]$. By
Lemma~\ref{lem:lifting}\textup{(iii)}, $z=\tilde x\tilde y$, and by
Lemma~\ref{lem:lifting}\textup{(i)}, the returned element is $xy$.

For the cost, forming $\hat x$ and $\hat y$ uses no arithmetic, since
each coefficient of $\tilde x$ and $\tilde y$ is already a polynomial of
degree less than $2k-1$ with coefficients in $\{0,1\}$. The product in
$E[L]$ uses the operations stated in
Proposition~\ref{prop:loop-correctness}
or~\ref{prop:central-correctness}, and $z$ has $|L|=2D$ coefficients,
each of degree at most $2k-2$ by Lemma~\ref{lem:lifting}\textup{(ii)}
and each reduced once modulo $2$ and $f$.
\end{proof}

Algorithm~\ref{alg:lift} is not an algorithm over $K$ in the sense of
Theorems~\ref{thm:main-a} and~\ref{thm:main-b}. Its arithmetic is
performed in $E$, whose elements have coefficients modulo a prime
exceeding $2kD$ and hence of more than $\log_2(2kD)$ bits, and no
bound on the number of arithmetic operations in $K$ required to
multiply in $K[L]$ is claimed here. Whether multiplication in $K[L]$ admits a
quasilinear algorithm in that model is not addressed in this paper.

\section{Implementation and benchmarks}
\label{sec:implementation}

\subsection{The implementation}
\label{subsec:loop-implementation}

The algorithms of this paper are implemented in \texttt{pycdloop}, a Python
package for arithmetic in the loops $Q_n$, their central products, and the loop
algebras of both over $\R$, $\C$, and the finite fields $\mathbb F_{p^k}$ with
$p\le2^{31}-1$ prime and $k\le12$. The package computes its Cayley--Dickson
products through \texttt{pyfastcd}, which accepts NumPy arrays
\cite{HarrisEtAl2020} and calls the \texttt{fastCD} library
\cite{LemleyFastCD2026}. The package is available on PyPI at
\href{https://pypi.org/project/pycdloop}{https://pypi.org/project/pycdloop},
with source code, tests, benchmark programs, reproduction scripts, and
archived data distributed under the MIT license at
\href{https://git.tanuki-cd.com/algebraity/pycdloop}{https://git.tanuki-cd.com/algebraity/pycdloop}
and mirrored on GitHub at
\href{https://github.com/TANUKI-CD/pycdloop}{https://github.com/TANUKI-CD/pycdloop}
\cite{LemleyPycdloop2026}. The
measurements reported in this section were taken on
November~15, 2026, using \texttt{pycdloop}~v0.1.0 and
\texttt{pyfastcd}~v1.4.0, and the recorded hashes of the measured
\texttt{pycdloop} source files match revision
\href{https://git.tanuki-cd.com/algebraity/pycdloop/commit/634c6bbf4d1cdcf73f858f6490ecb29515aa29c3}{\texttt{634c6bbf4d}},
which archives the
\href{https://github.com/TANUKI-CD/pycdloop/tree/634c6bbf4d1cdcf73f858f6490ecb29515aa29c3/benchmark-results/t440-20261115}{timing observations and reports} and the
\href{https://github.com/TANUKI-CD/pycdloop/tree/634c6bbf4d1cdcf73f858f6490ecb29515aa29c3/benchmark-results/stability-20261115}{accuracy summaries}.

The package provides more than one way to compute each product, and the
benchmarks of this section serve both to measure what
Algorithms~\ref{alg:loop} and~\ref{alg:central} gain over evaluating the
definition and, in the two settings where no single method is fastest
throughout, namely characteristic $2$ and central products, to evaluate the
defaults of the package. Throughout this section, rank refers to the rank of
the loop algebra, which is $2N$ for $R[Q_n]$ and $2D$ for $R[Q_n*Q_m]$.

Over fields of characteristic different from $2$, \texttt{split} uses the
decomposition and transforms of Algorithm~\ref{alg:loop} for $Q_n$ and
Algorithm~\ref{alg:central} for central products. Over $\R$, the quasilinear
kernel of \texttt{pyfastcd} uses the octonion base case of
\cite[Appendix~A]{Lemley2026} for $n\ge4$, which requires fewer operations than
the schedule for multiplication in $A_n$ used in this paper. The default
Cayley--Dickson kernel is quasilinear at every dimension $N\ge32$ over $\R$ and
$N\ge8$ over $\C$ and finite fields of odd characteristic
\cite{LemleyFastCD2026}. At smaller dimensions the default depends on the field,
the dimension, and the calling mode. At $N=16$ over $\R$, for example, the
measured build uses direct multiplication for single products and the algorithm
of Cariow and Cariowa \cite{CariowCariowa2015} for batches. The operation counts
in Propositions~\ref{prop:loop-correctness} and~\ref{prop:central-correctness}
use the level-two schedules $\mathcal E_n$ and, for the central product,
$\mathcal E_m$, while the timings below use the implementation defaults unless a
method is specified.

For $Q_n*Q_m$, the implementation takes the smaller level as the outer
one using the factor interchange of Section~\ref{sec:central}, so
$n\le m$. It computes the product in $A_n\otimes_R A_m$ by one of
three methods. The method \texttt{schedule} uses
Algorithm~\ref{alg:signed} with outer schedule $\mathcal E_n$ and
computes its products in $A_m$ through \texttt{pyfastcd}. It batches
sibling calls at each level of the outer recursion, so its auxiliary
storage exceeds the $O(D)$ bound of
Proposition~\ref{prop:central-correctness}, which refers to the
sequential evaluation. The method \texttt{rows} forms all
$N_n^2$ ordered products $U_pV_r$ in the formula preceding
Lemma~\ref{lem:coefficient-extension}. When $n=1$, \texttt{native}
computes one product in $A_m$ over $\C$ for base field $\R$, or over
$\mathbb F_{p^2}$ for base field $\mathbb F_p$ with $p\equiv3\pmod4$.

Over $\mathbb F_{2^k}$, the method \texttt{lift} uses
Algorithm~\ref{alg:lift} with $p=2^{31}-1$ on blocks of width at most
\[
h:=\begin{cases}
k & \text{if }1\le k\le6\,,\\
\lceil k/2\rceil & \text{if }7\le k\le12\,.
\end{cases}
\]
It uses an irreducible polynomial $\varphi\in\mathbb F_p[t]$ of degree
$2h-1$, so the auxiliary ring $E:=\mathbb F_p[t]/(\varphi)$ is a
field of that degree over $\mathbb F_p$. The largest extension degree
supported by \texttt{pyfastcd} is $12$. For a loop $L$, one product in
$E[L]$ suffices when $k\le6$. For $7\le k\le12$, the method splits
each coefficient into two blocks and combines three products in $E[L]$
by Karatsuba--Ofman multiplication \cite{KaratsubaOfman1962}.
For the middle product it adds the two blocks of each factor over
$\mathbb F_2$ before lifting. Every lifted block thus has coefficients
in $\{0,1\}$, so $h|L|<p$ suffices for recovery by
Lemma~\ref{lem:lifting}\textup{(ii)--(iii)}.

In every setting, the baseline is the algorithm \texttt{table}, which evaluates
$[x][y]=[xy]$ in the loop basis using $4D^2$ coefficient multiplications, where
$D=N$ for $Q_n$ and $D=N_nN_m$ for $Q_n*Q_m$, and which applies over every field
because it uses only the definition of the product. Finally, over fields of
characteristic different from $2$, an element $x\in R[Q_n]$ may also be kept in
\emph{split form}, the pair $(\Had_N\kappa_n(x),\varepsilon_n(x))$, and likewise
for central products. Since both components are multiplicative, two elements in
split form are multiplied by the coordinatewise product of their first
components together with the product in $A_n$ or $A_n \otimes_R A_m$ of their
second, with neither a transform nor a change of coordinates.

The benchmarks for \texttt{pycdloop} presented in this paper
ran on one logical CPU of an Intel Xeon Gold 6148 processor
(nominal frequency $2.40$ GHz) under Linux~7.2.11 with Python~3.14.8 and
NumPy~2.5.4. All timings were taken through the public Python interface and
therefore include the cost of each call, unlike the benchmarks of
\cite[Section~5]{Lemley2026}, which time native C calls; to separate this
fixed cost from the cost of the arithmetic, multiplication was measured both
singly, with one call for each operand pair, and in batches, with one call
for the complete set of pairs. Operands have independent standard normal
coefficients over $\R$ and $\C$ and independent uniformly distributed
residues over finite fields, and each timing benchmark consists of $24$ blocks
with a rotated task
order. Tables~\ref{tab:loop-real}--\ref{tab:central} report arithmetic
means of the block times per product and geometric means of the paired
per-block ratios. The archived reports give two-sided $95\%$ Student-$t$
intervals for these timings and ratios, and for ratios, the intervals are
computed on the logarithmic scale and then exponentiated.

Every computed product was checked against a reference that evaluates the
loop product from its definition, by exact equality over finite fields and,
over $\R$ and $\C$, by requiring the computed and reference products $z$ and
$z_{\rm ref}$ to satisfy
\[
\|z-z_{\rm ref}\|_\infty
\le 256\,\mathsf u D\,\|x\|_2\|y\|_2\,,
\]
where $2D$ is the loop order, $\mathsf u:=2^{-53}$ is the binary64
unit roundoff, and $\|\cdot\|_\infty$ is the largest absolute
coefficient. All $25{,}760$ checks passed.

\subsection{Cayley--Dickson loop algebras}
\label{subsec:loop-bench}

Table~\ref{tab:loop-real} gives the mean time per product in $\R[Q_n]$ in
nanoseconds for $1\le n\le10$. The algorithm \texttt{split} has the lower
mean time at every tested rank for both single and batched products, and its
advantage grows with the rank, reaching a factor of about $1560$ for batches
at rank $2048$. The number of multiplications accounts for only part of this
factor: at that rank, where $N=1024$, \texttt{table} uses $4N^2$
multiplications of two input-dependent quantities, while \texttt{split} uses
$N$ coordinatewise products and $37{,}258$ products in the real
\texttt{fastCD} schedule \cite[Table~7]{Lemley2026}, a ratio of about $110$.
The remainder is, aside from noise, a property of the two implementations in
\texttt{pycdloop}, in which \texttt{table} evaluates the definition through
NumPy, and not of the algorithms. At every tested rank from $128$
onward, \texttt{table} takes longer per product in a batch than in
single-product calls. For single products the fixed cost of a
Python call dominates at small ranks, where a single \texttt{split} product
takes $30.4$--$31.2$ microseconds at ranks $4$ through $32$, and the
single-product ratios accordingly remain below $4$ through rank $128$.

The benchmark also times the product in $A_n$ alone through
\texttt{pyfastcd} on the same data, which isolates the cost that
Algorithm~\ref{alg:loop} adds to that product. For batches, this overhead is
proportionally largest at small ranks, where Algorithm~\ref{alg:loop} takes
$19.9$ times as long as the product in $A_n$ at rank $4$, and it falls to a
factor of $1.94$ at rank $2048$, so that at the largest tested rank a
product in $\R[Q_n]$ costs about twice the product in $A_n$ that it
contains. Operands already in split form avoid the transforms and the
changes of coordinates, and the corresponding factors are $8.70$ and $1.22$.

\begin{table}[tbp]
\centering
\caption{Mean multiplication times in $\R[Q_n]$ in nanoseconds per product
over $24$ blocks, rounded to the nearest nanosecond. Ratio denotes the
paired ratio of the time of \texttt{table} to that of \texttt{split}, and
${}/A_n$ the paired ratio of the time of \texttt{split} to that of the
product in $A_n$ alone.}
\label{tab:loop-real}
\small
\setlength{\tabcolsep}{4pt}
\begin{tabular}{@{}rr rrrr rrr@{}}
\toprule
& & \multicolumn{4}{c}{Batch products}
& \multicolumn{3}{c}{Single products} \\
\cmidrule(lr){3-6}\cmidrule(l){7-9}
$2N$ & Pairs & \texttt{split} & \texttt{table} & Ratio & ${}/A_n$
& \texttt{split} & \texttt{table} & Ratio \\
\midrule
4 & 256 & 189 & 588 & 3.12 & 19.9 & 30923 & 46734 & 1.51 \\
8 & 256 & 203 & 820 & 4.04 & 17.6 & 30454 & 46970 & 1.54 \\
16 & 256 & 212 & 1682 & 7.95 & 10.9 & 30360 & 49839 & 1.64 \\
32 & 256 & 305 & 8950 & 29.3 & 4.91 & 31216 & 55352 & 1.77 \\
64 & 256 & 553 & 48355 & 87.5 & 2.98 & 31809 & 65007 & 2.04 \\
128 & 256 & 1128 & 197239 & 175 & 2.52 & 33777 & 121367 & 3.59 \\
256 & 64 & 2769 & 727264 & 263 & 2.50 & 35023 & 355695 & 10.2 \\
512 & 64 & 5685 & 3730925 & 658 & 2.19 & 38023 & 1345210 & 35.4 \\
1024 & 16 & 13475 & 14807170 & 1100 & 2.07 & 46132 & 5036054 & 109 \\
2048 & 16 & 28205 & 43788816 & 1560 & 1.94 & 59140 & 33963383 & 575 \\
\bottomrule
\end{tabular}
\end{table}

Table~\ref{tab:loop-fields} gives the corresponding batch times over $\C$ and
selected finite fields. Outside characteristic $2$ the comparison is as over
$\R$: over $\C$ and every tested field of odd characteristic, \texttt{split} has
the lower mean time at every tested rank for both single and batched products.
In characteristic $2$, by contrast, the outcome depends on both the degree of
the field and the rank. Over the seven tested fields $\mathbb F_{2^k}$ with
$k=1,2,4,6,7,8,12$, \texttt{lift} has the lower batch time at rank $8$ only for
$k=1,2,4$, with \texttt{table} up to $3.40$ times faster for the remaining
degrees, whereas at rank $512$ \texttt{lift} is faster in every tested field, by
factors from $33.6$ to $149$ for batches and from $11.1$ to $71.6$ for single
products. The package selects between the two exact methods by thresholds on the
field degree, the rank, and the number of products in a call. The default has
the lower mean time in $108$ of the $112$ tested combinations of field, rank,
and calling mode. In the other four, the default is \texttt{lift} and
\texttt{table} has the lower mean time.

\begin{table}[tbp]
\centering
\caption{Mean multiplication times in $F[Q_n]$ in nanoseconds per batched
product over $24$ blocks, rounded to the nearest nanosecond, for
\texttt{split} over $\C$ and fields of odd characteristic and for
\texttt{lift} over fields of characteristic $2$, each followed by the
paired ratio of the time of \texttt{table} to it. Batches hold $64$, $64$,
$16$, and $4$ products at the four ranks. The characteristic-$2$ rows
use the separate benchmark run over binary fields.}
\label{tab:loop-fields}
\small
\setlength{\tabcolsep}{4pt}
\begin{tabular}{@{}l rr rr rr rr@{}}
\toprule
& \multicolumn{2}{c}{$2N=8$}
& \multicolumn{2}{c}{$2N=32$}
& \multicolumn{2}{c}{$2N=128$}
& \multicolumn{2}{c}{$2N=512$} \\
\cmidrule(lr){2-3}\cmidrule(lr){4-5}\cmidrule(lr){6-7}\cmidrule(l){8-9}
Field & Time & Ratio & Time & Ratio & Time & Ratio & Time & Ratio \\
\midrule
$\C$ & 723 & 2.12 & 1215 & 6.17 & 5178 & 59.3 & 23892 & 182 \\
$\mathbb F_3$ & 1096 & 2.69 & 2151 & 10.0 & 8953 & 37.1 & 31247 & 177 \\
$\mathbb F_{2^{31}-1}$ & 1093 & 2.66 & 2244 & 8.65 & 8491 & 34.7 & 33122 & 150 \\
$\mathbb F_{3^5}$ & 3564 & 3.66 & 8285 & 17.1 & 31201 & 71.9 & 124150 & 324 \\
\midrule
$\mathbb F_2$ & 1226 & 1.07 & 2333 & 1.66 & 9159 & 6.11 & 37617 & 51.3 \\
$\mathbb F_{2^4}$ & 5869 & 1.30 & 18227 & 3.87 & 64712 & 20.0 & 276983 & 82.5 \\
$\mathbb F_{2^8}$ & 18026 & 0.457 & 55812 & 1.62 & 200685 & 7.05 & 792053 & 41.6 \\
\bottomrule
\end{tabular}
\end{table}

\FloatBarrier

\subsection{Central products}
\label{subsec:loop-bench-central}

Table~\ref{tab:central} gives the batch results for $Q_n*Q_m$ over $\R$.
Its first three timing columns give the product of two operands already in
split form, each including the $D$ coordinatewise products, so that they
differ only in the method used to multiply the component in
$A_n\otimes_\R A_m$; the last two give the time of
Algorithm~\ref{alg:central} on loop-basis operands with the default method
and the paired ratio of the time of \texttt{table} to it.

Where it applies, \texttt{native} has the lowest mean time. The comparison
between the other two methods follows the number of products in $A_m$, which is
$N_n^2$ for \texttt{rows} and $p_n$ for \texttt{schedule}. At $n=4$, these
numbers are $256$ and $178$, a saving too small to outweigh the linear
operations of \texttt{schedule} consistently, and for $2\le n\le4$ the two
methods are within a factor of $1.4$ of each other, with neither uniformly
faster. At $n=7$ the numbers are $16384$ and $3082$, and for $n\ge5$,
\texttt{schedule} has the lower mean time at every tested pair, with a paired
speedup of $6.07$ at $(7,7)$. The default is \texttt{native} where it applies,
\texttt{rows} for the remaining pairs with $n\le4$, and \texttt{schedule} for
$n\ge5$. For split-form operands at $(2,4)$ and $(4,8)$, \texttt{schedule} has
the lower mean time even though the default is \texttt{rows}. This crossover is
a property of the implementation, in which the linear operations of
\texttt{schedule} are carried out by NumPy, and could shift under another. With
these defaults, \texttt{split} has the lower mean time on loop-basis operands at
all $11$ pairs timed against \texttt{table}, with paired speedups from $1.88$ at
$(2,2)$ to $293$ at $(3,7)$.

\begin{table}[tbp]
\centering
\caption{Mean times in nanoseconds per batched product in $\R[Q_n*Q_m]$
over $24$ blocks, rounded to the nearest nanosecond. Split-form products use the
indicated method for the tensor component. Loop-basis products use
\texttt{split} with the default method. Ratio denotes the paired ratio of the
time of \texttt{table} to that of \texttt{split}. An unavailable method is
marked n/a and an unmeasured or unavailable entry by a dash.}
\label{tab:central}
\small
\setlength{\tabcolsep}{4pt}
\begin{tabular}{@{}lrr rrr rr@{}}
\toprule
& & & \multicolumn{3}{c}{Split-form operands}
& \multicolumn{2}{c}{Loop-basis operands} \\
\cmidrule(lr){4-6}\cmidrule(l){7-8}
$(n,m)$ & $2D$ & Pairs & \texttt{native} & \texttt{rows}
& \texttt{schedule} & \texttt{split} & Ratio \\
\midrule
$(1,4)$ & 64 & 64 & 1801 & 2323 & 2665 & 2239 & 13.9 \\
$(1,8)$ & 1024 & 4 & 34561 & 42656 & 51925 & 42746 & 243 \\
$(2,2)$ & 32 & 64 & n/a & 2434 & 2833 & 2809 & 1.88 \\
$(2,4)$ & 128 & 64 & n/a & 4763 & 4477 & 5458 & 21.7 \\
$(2,6)$ & 512 & 16 & n/a & 19402 & 20729 & 22726 & 126 \\
$(3,3)$ & 128 & 64 & n/a & 6362 & 8514 & 7084 & 15.8 \\
$(3,5)$ & 512 & 16 & n/a & 29071 & 37904 & 31974 & 88.6 \\
$(3,7)$ & 2048 & 4 & n/a & 131095 & 174554 & 150022 & 293 \\
$(4,4)$ & 512 & 16 & n/a & 52397 & 57376 & 55926 & 48.9 \\
$(4,6)$ & 2048 & 4 & n/a & 222998 & 271586 & 236249 & 184 \\
$(4,8)$ & 8192 & 4 & n/a & 911857 & 771936 & --- & --- \\
$(5,5)$ & 2048 & 4 & n/a & 382102 & 350873 & 361473 & 121 \\
$(5,7)$ & 8192 & 4 & n/a & 1680408 & 999042 & --- & --- \\
$(6,6)$ & 8192 & 4 & n/a & 3884691 & 1151597 & --- & --- \\
$(7,7)$ & 32768 & 2 & n/a & 34844261 & 5744267 & --- & --- \\
\bottomrule
\end{tabular}
\end{table}

\FloatBarrier

\subsection{Measured floating-point error}
\label{subsec:loop-stability}

The products computed over $\R$ and $\C$ are affected by rounding, and
Algorithm~\ref{alg:loop} rounds in its changes of coordinates and transforms
as well as in the product in $A_n$, so its accuracy is not evident from that
of \texttt{fastCD} alone. The accuracy of \texttt{split} and \texttt{table}
was therefore measured on $2{,}440$ operand pairs over $\R$ and $2{,}324$
over $\C$, drawn from $22$ real and $20$ complex families, against a
reference that evaluates every product exactly using integer arithmetic on
the stored binary64 inputs, as in the
\href{https://git.tanuki-cd.com/algebraity/fastCD/src/tag/v1.4.0/docs/stability.md}{\texttt{fastCD} accuracy study} \cite{LemleyFastCD2026}.
For nonzero $x,y$, the study reports
\[
\eta:=\frac{\|\widehat{xy}-xy\|_2}{\mathsf u\,\|x\|_2\|y\|_2}\,,
\]
where $\widehat{xy}$ is the computed product and the norms are Euclidean
norms of the coefficient vectors in the loop basis. As in
\cite[Section~5.5]{Lemley2026}, this quantity measures the error against the
operands, and the error relative to a product $xy\ne0$ is
\[
\mathsf u\,\eta\,\|x\|_2\|y\|_2/\|xy\|_2\,,
\]
which can be large when the
product is small relative to the operands.

Table~\ref{tab:loop-error} gives the median and the maximum of $\eta$ over the
pairs with Gaussian and uniform coordinates. Over $\R$, the median for
\texttt{split} rises from $0.78$ at rank $4$ to $2.73$ at rank $2048$, and a
least-squares fit of the squared medians for $1\le n\le10$ gives approximately
$0.81\,n-0.52$, so that over the tested range the median grows roughly as $\sqrt
n$. The error of \texttt{table} grows faster: it has the smaller median through
rank $64$, but from that rank onward its median grows by a factor of about $1.4$
per doubling, so that its square roughly doubles with the rank, and reaches
$9.68$ at rank $2048$. The difference is consistent with the way each algorithm
forms its sums. For each coefficient, \texttt{table} forms $N$ sums of two
products and uses NumPy to sum them along the non-contiguous axis indexed by
$p$, where rounding errors accumulate \cite{HighamMary2019}, whereas in
\texttt{split} the long sums are balanced: the Walsh--Hadamard butterflies form
a tree of depth $n$, and \texttt{fastCD} accumulates long sums in blocks whose
totals are combined by a balanced tree \cite{LemleyFastCD2026}. Over all
families and ranks, the largest value of $\eta$ observed for \texttt{split} is
$5.77$ over $\R$ and $4.29$ over $\C$. These figures describe only the tested
inputs for this particular implementation and thus cannot be used to make
general claims about accuracy, and no bound on the error of either algorithm is
proved in this paper.

\begin{table}[tbp]
\centering
\caption{Median and maximum of $\eta$ over the operand pairs with Gaussian
and uniform coordinates, for \texttt{split} and \texttt{table}. The pair
count applies to each domain measured at that rank. A dash marks a rank
omitted from the complex study.}
\label{tab:loop-error}
\small
\setlength{\tabcolsep}{4pt}
\begin{tabular}{@{}rr rr rr@{\hspace{12pt}}rr rr@{}}
\toprule
& & \multicolumn{4}{c}{Real coefficients}
& \multicolumn{4}{c}{Complex coefficients} \\
\cmidrule(lr){3-6}\cmidrule(l){7-10}
& & \multicolumn{2}{c}{\texttt{split}} & \multicolumn{2}{c}{\texttt{table}}
& \multicolumn{2}{c}{\texttt{split}} & \multicolumn{2}{c}{\texttt{table}} \\
\cmidrule(lr){3-4}\cmidrule(lr){5-6}\cmidrule(lr){7-8}\cmidrule(l){9-10}
$2N$ & Pairs & Median & Max & Median & Max & Median & Max & Median & Max \\
\midrule
4    & 256 & 0.78 & 2.81 & 0.67 & 1.50  & 0.99 & 2.84 & 0.79 & 1.53 \\
8    & 256 & 1.02 & 4.84 & 0.82 & 1.90  & 1.17 & 2.51 & 0.95 & 1.99 \\
16   & 256 & 1.28 & 2.74 & 1.04 & 2.00  & 1.61 & 3.26 & 1.14 & 1.86 \\
32   & 256 & 1.54 & 2.66 & 1.36 & 2.47  & 1.81 & 2.76 & 1.45 & 2.17 \\
64   & 256 & 1.90 & 3.05 & 1.86 & 2.76  & 2.02 & 2.63 & 1.91 & 2.45 \\
128  & 64  & 2.10 & 2.86 & 2.54 & 3.40  & 2.22 & 2.86 & 2.57 & 2.99 \\
256  & 64  & 2.25 & 2.61 & 3.58 & 4.37  & 2.38 & 2.63 & 3.59 & 3.93 \\
512  & 16  & 2.47 & 2.87 & 4.92 & 5.50  & 2.47 & 2.66 & 4.84 & 5.05 \\
1024 & 16  & 2.61 & 2.90 & 6.91 & 7.24  & \multicolumn{4}{c}{\textemdash} \\
2048 & 16  & 2.73 & 2.79 & 9.68 & 10.11 & \multicolumn{4}{c}{\textemdash} \\
\bottomrule
\end{tabular}
\end{table}

The splitting itself affects the accuracy in two ways. Because \texttt{split}
rounds $\kappa_n(x)$ and $\varepsilon_n(x)$ separately, an error in a small
component is measured against the size of the other, but on the four tested
families in which one component is small or zero, the median of $\eta$ is at
most $2.1$ times the dense-input median in Table~\ref{tab:loop-error} for the
same field and rank. Keeping a left-nested product in split form avoids rounding
in the intermediate conversions to and from the loop basis. For the $712$ real
chains of four or eight factors, the error with this evaluation is smaller than
with conversion to the loop basis after every product in $540$ cases, equal in
$18$, and larger in $154$. The corresponding counts for the $672$ complex chains
are $518$, $5$, and $149$.

\section{Outlook and applications}
\label{sec:outlook}

This paper gives quasilinear algorithms for multiplication in $R[Q_n]$ and
$R[Q_n*Q_m]$ and benchmarks their implementation in \texttt{pycdloop}. It
remains open whether these algorithms are optimal in a specified computational
model, and whether multiplication in $A_n\otimes_RA_m$, and hence in
$R[Q_n*Q_m]$, can be computed using $O(D\log D)$ operations when both $n$ and
$m$ grow. Over finite fields $K$ of characteristic $2$, it is likewise open
whether multiplication in $K[Q_n]$ and $K[Q_n*Q_m]$ admits a quasilinear
algorithm counted in operations over $K$. Furthermore, while loop algebras have
been applied in various computational settings, including cryptography and
coding theory \cite{GonzalezEtAl2004, MarkovEtAl2020}, potential applications
of Cayley--Dickson loop algebras do not appear to have been explored in detail;
it could be fruitful to consider possible applications of these objects in the
computer sciences and other technical fields which may benefit from their
algebraic structure. Finally, the algebras presented in this paper are not the
only algebras constructed from Cayley--Dickson algebras, and it would be
interesting to see to what other algebras the quasilinear algorithm for $A_n$
could be extended or applied.

\section{Statements and disclosures}
\label{sec:disclosures}

\subsection{AI use declaration}

Generative AI is a substantial part of the author's research process. The author
made use of frontier models through OpenAI's ChatGPT and Codex and Anthropic's
Claude and Claude Code in the development of this paper.

These models were used to assist with literature search and review, identify
open problems, suggest research directions, explain and analyze papers, assist
with drafting and editing, discuss proof ideas, suggest lemmas and consequences,
outline proofs, and develop computational tests of results. They were also used
substantively during the development of proofs, in some cases proposing key
arguments or generating candidate proofs in response to mathematical problems
posed by the author.

\subsection{Code availability}

The code and data discussed in this paper and used to support its results are
publicly available under the MIT license at
\href{https://git.tanuki-cd.com/algebraity/pycdloop}{https://git.tanuki-cd.com/algebraity/pycdloop},
with a GitHub mirror available at
\href{https://github.com/TANUKI-CD/pycdloop}{https://github.com/TANUKI-CD/pycdloop}.

\subsection{Competing interests and funding}

The author is the developer of the \texttt{pycdloop} package benchmarked in
Section~\ref{sec:implementation} and declares no other competing interests.
This research received no specific funding.

\subsection{Statement of responsibility}

All results presented in this paper have been independently verified by the
author. The author takes full responsibility for the contents of the paper,
including its mathematical claims, arguments, citations, and reported
computations.

\printbibliography

\end{document}